\documentclass[11pt]{article}
\usepackage{amssymb}
\usepackage{times}

\title{A dichotomy characterizing analytic digraphs of uncountable Borel chromatic number in any dimension.\indent}
\author{Dominique LECOMTE}
\date{\it ~Trans. Amer. Math. Soc.\rm ~361 (2009), 4181-4193}

\def\ufootnote#1{\let\savedthfn\thefootnote\let\thefootnote\relax
\footnote{#1}\let\thefootnote\savedthfn\addtocounter{footnote}{-1}}

\newcommand{\Ana}{{\it\Sigma}^{1}_{1}}
\newcommand{\Ca}{{\it\Pi}^{1}_{1}}
\newcommand{\Borel}{{\it\Delta}^{1}_{1}}

\newcommand{\ana}{{\bf\Sigma}^{1}_{1}}

\newtheorem{thm} {Theorem} [section]
\newtheorem{defi} [thm] {Definition}

\newtheorem{lem} [thm] {Lemma}

\begin{document}

\maketitle

\ufootnote{{\it 2000 Mathematics Subject Classification.}~Primary: 03E15, Secondary: 54H05}

\ufootnote{{\it Keywords and phrases.}~Borel chromatic number, dimension}

\ufootnote{{\bf Acknowledgements.}~I would like to thank the anonymous referee for the simplification of some proofs in this paper.}

\noindent {\bf Abstract.} We study the extension of the Kechris-Solecki-Todor\v cevi\'c dichotomy on analytic graphs to dimensions higher than 2. We prove that the extension is possible in any dimension, finite or infinite. The original proof works in the case of the finite dimension. We first prove that the natural extension does not work in the case of the infinite dimension, for the notion of continuous homomorphism used in the original theorem. Then we solve the problem in the case of the infinite dimension. Finally, we prove that the natural extension works in the case of the infinite dimension, but for the notion of Baire-measurable homomorphism.

\section{$\!\!\!\!\!\!$ Introduction.}\indent

 The reader should see [K] for the standard descriptive set theoretic notation used in this paper. We study a definable coloring problem, in any dimension. We will need some more notation:\bigskip
 
\noindent\bf Notation.\rm ~In this paper, $2\!\leq\! d\!\leq\!\omega$ will be a cardinal, i.e., any dimension of an actual product making sense in the context of descriptive set theory. The letters $X$, $Y$ will refer to some sets. We set 
$$\Delta^d(X)\! :=\!\{ (x_i)_{i\in d}\!\in\! X^d\mid\forall i\!\in\! d\ \ x_i\! =\! x_0\}.$$

\begin{defi} Let $A\!\subseteq\! X^d$. We say that $A$ is a $\underline{digraph}$ if 
$A\cap\Delta^d(X)\! =\!\emptyset$.\end{defi}

\noindent\bf Notation.\rm ~Let $u\! :\! X\!\rightarrow\! Y$ be a map. We define a map $u^d\! :\! X^d\!\rightarrow\! Y^d$ by 
$$u^d[(x_i)_{i\in d}]\! :=\! [u(x_i)]_{i\in d}.$$

\begin{defi} Let $A\!\subseteq\! X^d$ be a digraph.\smallskip

\noindent (a) A $\underline{coloring}$ of $[X,A]$ is a map 
$c\! :\! X\!\rightarrow Y$ such that $A\cap (c^d)^{-1}[\Delta^d(Y)]\! =\!\emptyset$.\smallskip

\noindent (b) Assume that $X$ is a Polish space. The $\underline{Borel\ chromatic\ number}$ of $[X,A]$ is 
$$\chi_B(A)\! :=\!\hbox{\it min}\{\ \hbox{\it Card}(Y)\mid Y\ \hbox{\it is\ a\ Polish\ space\ and\ there\ is\ a\ Borel\ coloring}\ c\! :\! X\!\rightarrow Y\ \hbox{\it of}\ [X,A]\ \}.$$\end{defi}

\vfill\eject

 The goal of this paper is to characterize the analytic digraphs of uncountable Borel chromatic 
number. This has been done in [K-S-T] for graphs, i.e., for symmetric digraphs, when 
$d\! =\! 2$. We will give such a characterization in terms of the following notion of comparison between relations.\bigskip

\noindent\bf Notation.\rm ~Assume that $X$, $Y$ are Polish spaces, and let $A$ (resp., $B$) be a subset of $X^d$ (resp., $Y^d$). We set
$$[X,A]\preceq_{B}[Y,B]~\Leftrightarrow ~\exists u\! :\! X\!\rightarrow\! Y\ \ \hbox{Borel\ with}\ 
A\!\subseteq\! (u^d)^{-1}(B).$$
In this case, we say that $u$ is a Borel $\underline{homomorphism}$ from $[X,A]$ into $[Y,B]$. This notion essentially makes sense for digraphs (we can take $u$ to be constant if $B$ is not a digraph). If $u$ is continuous (resp., Baire-measurable, arbitrary), then we write $\preceq_{c}$ (resp., $\preceq_{Bm}$, $\preceq$) instead of $\preceq_{B}$. Note that $\chi_B(A)\!\leq\!\omega$ is equivalent to $[X,A]\preceq_B[\omega ,\neg\Delta^d(\omega )]$.\bigskip

 We also have to introduce minimum digraphs of uncountable Borel chromatic number:\bigskip

\noindent $\bullet$ Let $\psi_d\! :\!\omega\!\rightarrow\! d^{<\omega}$ be the natural bijection, for 
$d\!\leq\!\omega$. More specifically,\bigskip

- If $d\! <\!\omega$, then $\psi_d (0)\! :=\!\emptyset$ is the sequence of length $0$, $\psi_d (1)\! :=\! 0$, ..., $\psi_d (d)\! :=\! d\! -\! 1$ are the sequences of length $1$, and so on.\bigskip 

- If $d\! =\!\omega$, then let $(p_n)_{n\in\omega}$ be the sequence of prime numbers, and 
$I\! :\!\omega^{<\omega}\!\rightarrow\!\omega$ defined by $I(\emptyset )\! :=\! 1$, and 
$I(s)\! :=\! p_0^{s(0)+1}...p_{\vert s\vert -1}^{s(\vert s\vert -1)+1}$ if $s\!\not=\!\emptyset$. Note that $I$ is one-to-one, so that there is an increasing bijection $\varphi\! :\!\hbox{\rm Seq}\! :=\! I[\omega^{<\omega}]\!\rightarrow\!\omega$. If $t\!\in\!\hbox{\rm Seq}\!\subseteq\!\omega$, then we will denote 
$\overline{t}\! :=\! I^{-1}(t)\!\in\!\omega^{<\omega}$. We set $\psi_\omega\! :=\! (\varphi\circ I)^{-1}\! :\!\omega\!\rightarrow\!\omega^{<\omega}$. Note that $\psi_\omega$ is a bijection.\bigskip

\noindent $\bullet$ Note also that $|\psi_d (n)|\!\leq\! n$ if $n\!\in\!\omega$. Indeed, this is clear if $d\! <\!\omega$. If $d\! =\!\omega$, then
$$I[\psi_\omega (n)\vert 0]\! <\! I[\psi_\omega (n)\vert 1]\! <\! ...\! <\! 
I[\psi_\omega (n)]\hbox{\rm ,}$$
so that $(\varphi\circ I)[\psi_\omega (n)\vert 0]\! <\! (\varphi\circ I)[\psi_\omega (n)\vert 1]\! <\! ...\! <\! (\varphi\circ I)[\psi_\omega (n)]\! =\! n$. This implies that $|\psi_\omega (n)|\!\leq\! n$.\bigskip

\noindent $\bullet$ Let $n\!\in\!\omega$. As $|\psi_d (n)|\!\leq\! n$, we can define 
$s^d_{n}\! :=\!\psi_d (n)0^{n-|\psi_d (n)|}$. The crucial properties of the sequence $(s^d_{n})_{n\in\omega}$ are the following:\bigskip

- For each $s\!\in\! d^{<\omega}$, there is $n\!\in\!\omega$ such that $s\!\subseteq\! s^d_{n}$ (we say that 
$(s^d_{n})_{n\in\omega}$ is $\underline{dense}$ in $d^{<\omega}$).\bigskip

- $|s^d_{n}|\! =\! n$.\bigskip

\noindent $\bullet$ We put
$$\mathbb{A}_{d}\! :=\!\{ (s^d_{n}i\gamma )_{i\in d}\mid n\!\in\!{\omega}\ \hbox{\rm and}\ \gamma\!\in\! d^{\omega}\}\!\subseteq\! (d^\omega)^d.$$
\noindent Note that $\mathbb{A}_{d}\!\in\!\ana$ since the map 
$(n,\gamma )\!\mapsto\! (s^d_{n}i\gamma )_{i\in d}$ is continuous.\bigskip

 The previous definitions were given, when $d\! =\! 2$, in [K-S-T], where the following is proved:

\begin{thm} (Kechris, Solecki, Todor\v cevi\'c) Let $X$ be a Polish space, and $A\!\in\!\ana (X^2)$. Then exactly one of the following holds:\smallskip  

\noindent (a) $[X,A]\preceq_{B}[{\omega},\neg\Delta^2(\omega )]$.\smallskip  

\noindent (b) $[2^{\omega},\mathbb{A}_{2}]\preceq_{c}[X,A]$.\end{thm}
 
 This result can be extended to any finite dimension $d$, with the same proof.
 
\vfill\eject
 
\begin{thm} Let $d\!\geq\! 2$ be an integer, $X$ a Polish space, and $A\!\in\!\ana (X^d)$. Then exactly one of the following holds:\smallskip  

\noindent (a) $[X,A]\preceq_{B}[{\omega},\neg\Delta^d(\omega )]$.\smallskip  

\noindent (b) $[d^{\omega},\mathbb{A}_{d}]\preceq_{c}[X,A]$.\end{thm}

 We want to study the case of the infinite dimension. 
 
\begin{thm} We cannot extend Theorem 1.4 to the case where $d\! =\!\omega$.\end{thm}
 
\noindent\bf Notation.\rm ~In order to get a positive result in the case of the infinite dimension, we put
$$\mathbb{G}\! :=\!\{\alpha\!\in\!\omega^\omega\mid\forall m\!\in\!\omega\ \exists n\!\geq\! m\ \ s^\omega_n0\!\subseteq\!\alpha\}.$$
Note that $\mathbb{G}$ is a dense $G_\delta$ subset of $\omega^\omega$.\bigskip

 The main result of this paper is the following:
   
\begin{thm} Let $X$ be a Polish space, and $A\!\in\!\ana (X^\omega )$. Then exactly one of the following holds:\smallskip  

\noindent (a) $[X,A]\preceq_{B}[{\omega},\neg\Delta^\omega (\omega )]$.\smallskip  

\noindent (b) $[\mathbb{G},\mathbb{A}_{\omega}\cap\mathbb{G}^\omega ]\preceq_{c}[X,A]$.\end{thm}

 So we have a general characterization, in any dimension $d$, of analytic relations $A\!\subseteq\! X^d$ for which $[X,A]\not\preceq_{B}[{\omega},\neg\Delta^d(\omega )]$. In particular, we have a 
characterization of analytic digraphs of uncountable Borel chromatic number.\bigskip

 Theorem 1.5 says that we cannot extend Theorem 1.4 to the case where $d\! =\!\omega$ for the notion of continuous homomorphism in (b). However, the extension of Theorem 1.4 to the case where 
 $d\! =\!\omega$ is possible for the notion of Baire-measurable homomorphism:

\begin{thm} Let $X$ be a Polish space, and $A\!\in\!\ana (X^\omega )$. Then exactly one of the following holds:\smallskip  

\noindent (a) $[X,A]\preceq_{B}[{\omega},\neg\Delta^\omega (\omega )]$.\smallskip  

\noindent (b) $[\omega^{\omega},\mathbb{A}_{\omega}]\preceq_{Bm}[X,A]$.\end{thm}

\section{$\!\!\!\!\!\!$ The proof in finite dimension.}\indent

 Let us start with two general lemmas:

\begin{lem} Let $G$ be a dense $G_\delta$ subset of $d^\omega$. Then 
$[G,\mathbb{A}_{d}\cap G^d]\not\preceq_{Bm}[{\omega},\neg\Delta^d(\omega )]$.\end{lem}

\noindent\bf Proof.\rm\ We argue by contradiction. This gives a Baire-measurable function 
$u\! :\! G\!\rightarrow\!\omega$ such that 
$\mathbb{A}_d\cap G^d\!\subseteq\! (u^d)^{-1}[\neg\Delta^d(\omega)]$. As 
$G\! =\!\bigcup_{i\in\omega}\ u^{-1}(\{ i\})$, there is an integer $i_0$ such that 
$u^{-1}(\{ i_0\})$ is not meager and has the Baire property in $G$. This implies the existence of $s\!\in\! d^{<\omega}$ such that $(G\cap N_s)\!\setminus\! u^{-1}(\{ i_0\})$ is meager. Let $H$ be a dense $G_\delta$ subset of $G$ such that $H\cap N_{s}\subseteq\! u^{-1}(\{ i_0\})$. We choose $n\!\in\!\omega$ with 
$s\subseteq s^d_{n}$. Note that $f^i_n\! :\! N_{s^d_{n}0}\!\rightarrow\! N_{s^d_{n}i}$ defined by 
$f^i_n (s^d_{n}0\gamma )\! :=\! s^d_{n}i\gamma$ is an homeomorphism. This implies that 
$\bigcap_{i\in\omega}\ (f^i_n)^{-1}(H)$ is a dense $G_\delta$ subset of $N_{s^d_{n}0}$. We choose $s^d_{n}0\gamma\!\in\!\bigcap_{i\in\omega}\ (f^i_n)^{-1}(H)$. We get 
$(s^d_{n}i\gamma )_{i\in d}\!\in\!\mathbb{A}_d\cap (H\cap N_{s})^d\!\subseteq\! [u^{-1}(\{ i_0\})]^d$, which contradicts the fact that $\mathbb{A}_d\cap G^d\!\subseteq\! (u^d)^{-1}[\neg\Delta^d(\omega)]$.
\hfill{$\square$}

\vfill\eject

\begin{defi} Let $A\!\subseteq\! X^d$. We say that $C\!\subseteq\! X$ is $\underline{A\! -\! discrete}$ if $A\cap C^d\! =\!\emptyset$.\end{defi}

\noindent\bf Notation.\rm\ The reader should see [M] for the basic notions of effective descriptive set theory. Assume that $X$ and $X^d$ are recursively presented Polish spaces, and that 
$A\!\in\!\Ana (X^d)$. We put 
$$U\! :=\!\bigcup\ \{ D\!\in\!\Borel (X)\mid D\ \hbox{is}\ A\hbox{\rm -discrete}\}.$$
Note that $U\!\in\!\Ca (X)$ if the projections are recursive.

\begin{lem} Assume that $X$ and $X^d$ are recursively presented Polish spaces, 
$A\!\in\!\Ana (X^d)$, and $U\! =\! X$. Then 
$[X,A]\preceq_{B}[{\omega},\neg\Delta^d(\omega )]$.\end{lem}

\noindent\bf Proof.\rm\ As $U\! =\! X$, there is a partition $(D_n)_{n\in\omega}$ of $X$ into $A$-discrete 
$\Borel$ sets. We define a Borel map $u\! :\! X\!\rightarrow\!\omega$ by 
$u(x)\! =\! n\Leftrightarrow x\!\in\! D_n$. If $(x_i)_{i\in d}\!\in\! A$, then we cannot have 
$[u(x_i)]_{i\in d}\!\in\!\Delta^d(\omega )$, since the $D_n$'s are $A$-discrete.\hfill{$\square$}\bigskip

 We will recall the proof of Theorem 1.4, to show the problem appearing in the case of the infinite dimension. It is essentially identical to the one in [K-S-T], except that we do not use Choquet games.\bigskip
 
\noindent\bf Notation.\rm\ Let $Z$ be a recursively presented Polish space. The 
$\underline{Gandy\! -\! Harrington\ topology}$ on $Z$ is generated by $\Ana (Z)$ and denoted ${\it\Sigma}_Z$. It is finer than the initial topology of $Z$, so that $[Z ,{\it\Sigma}_Z ]$ is $T_1$. 
As $\Ana (Z)$ is countable (see 3F.6 in [M]), $[Z ,{\it\Sigma}_Z ]$ is second countable. We set  
$$\Omega_Z :=\{ z\!\in\! Z\mid\omega_1^z\! =\!\omega_1^{\hbox{\rm CK}}\}.$$ 
Recall that $\Omega_Z$ is $\Ana (Z)$ and dense in $[Z,{\it\Sigma}_Z]$ (see III.1.5 in [S]; the second assertion is Gandy's basis theorem). Recall also that $W\cap\Omega_Z$ is a clopen subset of $[\Omega_Z ,{\it\Sigma}_Z ]$ for each $W\!\in\!\Ana (Z)$. Indeed, it is obviously open. Let 
$f\! :\! Z\!\rightarrow\!\omega^\omega$ be $\Borel$ such that 
$Z\!\setminus\! (W\cap\Omega_Z)\! =\! f^{-1}(WO)$ (see 4A.3 in [M]). We get
$$z\!\in\!\Omega_Z\!\setminus\! (W\cap\Omega_Z)\ \Leftrightarrow\ z\!\in\!\Omega_Z\ \ \hbox{\rm and}\ \ 
\exists\xi\! <\!\omega_1^{\hbox{\rm CK}}\ (f(z)\!\in\! WO\ \ \hbox{\rm and}\ \ \vert f(z)\vert\!\leq\!\xi ).$$
This proves that $W\cap\Omega_Z$ is closed (see 4A.2 in [M]). In particular, 
$[\Omega_Z ,{\it\Sigma}_Z ]$ is zero-dimensional, and regular. By Theorem 4.2 in [H-K-L] and 8.16.(iii) in [K], $[\Omega_Z ,{\it\Sigma}_Z ]$ is strong Choquet. By 8.18 in [K], $[\Omega_Z ,{\it\Sigma}_Z ]$ is a Polish space. So we fix a complete compatible metric $d_Z$ on $[\Omega_{Z},{\it\Sigma}_{Z}]$.\bigskip

\noindent\bf Proof of Theorem 1.4.\rm\ Note first that we cannot have (a) and (b) simultaneously, by Lemma 2.1.\bigskip

\noindent $\bullet$ We may assume that $X$ is a recursively presented Polish space and that 
$A\!\in\!\Ana (X^d)$. We set $\Phi\! :=\!\{ C\!\subseteq\! X\mid C\ \hbox{\rm is}\ A\hbox{\rm -discrete}\}$. As 
$\Phi$ is $\Ca$ on $\Ana$, the first reflection theorem ensures that if $C\!\in\!\Ana (X)$ is $A$-discrete, then there is $D\!\in\!\Borel (X)$ which is $A$-discrete and contains $C$ (see 35.C in [K]).\bigskip

\noindent $\bullet$ By Lemma 2.3 we may assume that $U\!\not=\! X$, so that $Y\! :=\! X\!\setminus\! U$ is a nonempty $\Ana$ subset of $X$. The previous point gives the following key property:
$$\forall C\!\in\!\Ana (X)\ \ \ (\emptyset\!\not=\! C\!\subseteq\! Y\ \Rightarrow\ A\cap C^d\!\not=\!
\emptyset ).$$

\vfill\eject

\noindent $\bullet$ We construct $(x_s)_{s\in d^{<\omega}}\!\subseteq\! Y$, 
$(V_s)_{s\in d^{<\omega}}\!\subseteq\!\Ana (X)$ and 
$(U_{n,t})_{(n,t)\in\omega\times d^{<\omega}}\!\subseteq\!\Ana (X^d)$ satisfying the following conditions:
$$\begin{array}{ll}
& (1)\ x_s\!\in\! V_s\!\subseteq\! Y\cap\Omega_X\ \ \hbox{\rm and}\ \ 
(x_{s_n^dit})_{i\in d}\!\in\! U_{n,t}\!\subseteq\! A\cap Y^d\cap
\Omega_{X^d}\hbox{\rm ,}\cr & \cr
& (2)\ V_{sm}\!\subseteq\! V_{s}\ \ \hbox{\rm and}\ \ 
U_{n,tm}\!\subseteq\! U_{n,t}\hbox{\rm ,}\cr & \cr
& (3)\ \hbox{\rm diam}_{d_X}(V_s)\!\leq\! 2^{-\vert s\vert}\ \ \hbox{\rm and}\ \ 
\hbox{\rm diam}_{d_{X^d}}(U_{n,t})\!\leq\! 2^{-n-1-\vert t\vert}.
\end{array}$$
$\bullet$ Assume that this is done. Fix $\alpha\!\in\! d^\omega$. Then 
$(V_{\alpha\vert p})_{p\in\omega}$ is a decreasing sequence of nonempty clopen subsets of 
$[\Omega_X ,{\it\Sigma}_X ]$ whose $d_X$-diameters tend to zero, so there is $u(\alpha )$ in their intersection. This defines $u\! :\! d^\omega\!\rightarrow\! X$. Note that 
$d_X[x_{\alpha\vert p},u(\alpha )]\!\leq\!\hbox{\rm diam}_{d_X}(V_{\alpha\vert p})\!\leq\! 2^{-p}$, so that $u$ is continuous and $(x_{\alpha\vert p})_{p\in\omega}$ tends to $u(\alpha )$ in $[X,{\it\Sigma}_X]$.\bigskip

 If $(s^d_ni\gamma )_{i\in d}\!\in\!\mathbb{A}_d$, then 
$(U_{n,\gamma\vert p})_{p\in\omega}$ is a decreasing sequence of nonempty clopen subsets of $[\Omega_{X^d} ,{\it\Sigma}_{X^d}]$ whose $d_{X^d}$-diameters tend to zero, so there is 
$(\alpha_i)_{i\in d}$ in their intersection. Note that $(\alpha_i)_{i\in d}\!\in\! A$. Moreover, the sequence 
$([x_{s^d_ni(\gamma\vert p)}]_{i\in d})_{p\in\omega}$ tends to $(\alpha_i)_{i\in d}$ in 
$[X^d,{\it\Sigma}_{X^d}]$, and in $[X^d,{\it\Sigma}_X^d]$ too. As 
$(x_{s^d_ni(\gamma\vert p)})_{p\in\omega}$ tends to $u(s^d_ni\gamma )$, we get 
$u(s^d_ni\gamma )\! =\!\alpha_i$, for each $i\!\in\! d$. Thus $[u(s^d_ni\gamma )]_{i\in d}\!\in\! A$.\bigskip

\noindent $\bullet$ So it is enough to see that the construction is possible. As $Y$ is a nonempty $\Ana$ subset of $X$, we can choose $x_\emptyset\!\in\! Y\cap\Omega_X$, and $V_\emptyset\!\in\!\Ana (X)$ such that $x_\emptyset\!\in\! V_\emptyset\!\subseteq\! Y\cap\Omega_X$ and 
$\hbox{\rm diam}_{d_X}(V_\emptyset )\!\leq\! 1$. Assume that $(x_s)_{\vert s\vert\leq l}$, 
$(V_s)_{\vert s\vert\leq l}$ and $(U_{n,t})_{n+1+\vert t\vert\leq l}$ satisfying (1)-(3) have been constructed, which is the case for $l\! =\! 0$. Let $C$ be the following set:
$$\{ x\!\in\! X\mid\exists (y_s)_{s\in d^l}\!\in\! X^{d^l}\ \ y_{s^d_l}\! =\! x\ \ 
\hbox{\rm and}\ \ \forall s\!\in\! d^l\ y_s\!\in\! V_s\ \ \hbox{\rm and}\ \ \forall n\! <\! l\ \ 
\forall t\!\in\! d^{l-n-1}\ (y_{s_n^dit})_{i\in d}\!\in\! U_{n,t}\}.$$
Then $C\!\in\!\Ana (X)$\bf\ since $d$ is an integer\rm , $x_{s^d_l}\!\in\! C\!\subseteq\! Y$ by induction assumption. So there is $(x_{s^d_li})_{i\in d}$ in $A\cap C^d\cap\Omega_{X^d}$, by the key property. As $x_{s^d_lm}\!\in\! C$, we get $(x_{sm})_{s\in d^l\setminus\{ s^d_l\}}$. It remains to  choose\bigskip 

- $V_{sm}\!\in\!\Ana (X)$ with $x_{sm}\!\in\! V_{sm}\!\subseteq\! V_s$ and 
$\hbox{\rm diam}_{d_X}(V_{sm})\!\leq\! 2^{-l-1}$, for $s\!\in\! d^l$ and $m\!\in\! d$.\bigskip

- $U_{l,\emptyset}\!\in\!\Ana (X^d)$ with  
$(x_{s^d_li})_{i\in d}\!\in\! U_{l,\emptyset}\!\subseteq\! A\cap Y^d\cap\Omega_{X^d}$ and 
$\hbox{\rm diam}_{d_{X^d}}(U_{l,\emptyset})\!\leq\! 2^{-l-1}$.\bigskip

- $U_{n,tm}\!\in\!\Ana (X^d)$ with  
$(x_{s_n^ditm})_{i\in d}\!\in\! U_{n,tm}\!\subseteq\! U_{n,t}$ and 
$\hbox{\rm diam}_{d_{X^d}}(U_{n,tm})\!\leq\! 2^{-l-1}$, for 
$(n,t)$ in $\omega\!\times\! d^{<\omega}$ with $n\! +\! 1\! +\!\vert t\vert\! =\! l$ and 
$m\!\in\! d$.\hfill{$\square$}

\section{$\!\!\!\!\!\!$ The natural extension in infinite dimension does not work.}\indent 

 Theorem 1.5 is a consequence of Lemma 2.1 and of the following result:\bigskip
  
\noindent\bf Theorem 3\it\ $[\omega^\omega,\mathbb{A}_{\omega}]\not\preceq_{c}
[\mathbb{G},\mathbb{A}_{\omega}\cap\mathbb{G}^\omega ]$.

\vfill\eject

\noindent\bf Proof.\rm\ We argue by contradiction. This gives a continuous map 
$u\! :\!\omega^\omega\!\rightarrow\!\mathbb{G}$ with 
$\mathbb{A}_{\omega}\!\subseteq\! (u^\omega )^{-1}(\mathbb{A}_{\omega})$.\bigskip

\noindent $\bullet$ Let us prove that there is $\alpha\!\in\!\omega^\omega$ and 
$(s_n)_{n\in\omega}\!\in\! (\omega^{<\omega})^\omega$ such that 
$$u[\beta (0)0^{\alpha (0)}\beta (1)0^{\alpha (1)}...]\! =\! s_0\beta (0)s_1\beta (1)...$$
for each $\beta\!\in\!\omega^\omega$. We construct $\alpha (n)$ and $s_n$ by induction on $n$. Assume that $\alpha\vert n$ and $(s_p)_{p<n}$ are constructed satisfying 
$$s^\omega_{\Sigma_{j\leq p} [1+\alpha (j)]}\!\subseteq\! 0^\infty\ \ \hbox{\rm and}\ \ 
[t(0)0^{\alpha (0)}...t(p)0^{\alpha (p)}\!\subseteq\!\gamma\ \Rightarrow\ s_0t(0)...s_pt(p)\!\subseteq\! u(\gamma )]$$
for each $p\! <\! n$ and $t\!\in\!\omega^{p+1}$. We will construct $\alpha (n)$ and $s_n$ satisfying 
$$s^\omega_{\Sigma_{j\leq n} [1+\alpha (j)]}\!\subseteq\! 0^\infty\ \ \hbox{\rm and}\ \ 
[t(0)0^{\alpha (0)}...t(n)0^{\alpha (n)}\!\subseteq\!\gamma\ \Rightarrow\ s_0t(0)...s_nt(n)\!\subseteq\! u(\gamma )]$$
for each $t\!\in\!\omega^{n+1}$, which will be enough. Note first that 
there are $m\!\in\!\omega$ and $\delta\!\in\!\omega^\omega$ with 
$[u(s^\omega_{\Sigma_{j<n} [1+\alpha (j)]}i0^\infty )]_{i\in\omega}\! =\! 
(s^\omega_mi\delta )_{i\in\omega}$. As $u$ is continuous, there is $p\!\in\!\omega$ such that
$$\begin{array}{ll}
s^\omega_{\Sigma_{j<n} [1+\alpha (j)]}0^{p+1}\!\!\!\! & \subseteq\!\gamma\ \Rightarrow\ 
s^\omega_m0\!\subseteq\! u(\gamma )\hbox{\rm ,}\cr & \cr
s^\omega_{\Sigma_{j<n} [1+\alpha (j)]}10^{p}\!\!\!\! & \subseteq\!\gamma\ \Rightarrow\ 
s^\omega_m1\!\subseteq\! u(\gamma ).
\end{array}$$
Note that $s^\omega_{\Sigma_{j<n} [1+\alpha (j)]}i0^{p}\!\subseteq\!\gamma\ \Rightarrow\ 
s^\omega_mi\!\subseteq\! u(\gamma )$, for each $i\!\in\!\omega$. Indeed, let 
$\varepsilon\!\in\!\omega^\omega$. Then 
$[u(s^\omega_{\Sigma_{j<n} [1+\alpha (j)]}i0^{p}\varepsilon )]_{i\in\omega}\!\in\!
\mathbb{A}_{\omega}\cap [N_{s^\omega_m0}\!\times\! N_{s^\omega_m1}\!\times\! (\omega^\omega )^\omega ]\!\subseteq\!\Pi_{i\in\omega}\ N_{s^\omega_mi}$. In particular, this implies that 
$s_00...s_{n-1}0\!\subseteq\! s^\omega_m$ since $s_00...s_{n-1}0\!\subseteq\! 
u(s^\omega_{\Sigma_{j<n} [1+\alpha (j)]}i0^{p}\varepsilon )$.\bigskip

\noindent - If $n\! =\! 0$, then we choose 
$\alpha (0)\!\geq\! p$ such that $0^{1+\alpha (0)}\! =\! s^\omega_{1+\alpha (0)}$, we set 
$s_0\! :=\! s^\omega_m$, and we are done.\bigskip

\noindent - If $n\! >\! 0$, then we set 
$s_n\! :=\! s^\omega_m\! -\! (s_00...s_{n-1}0)$. We will prove, by induction on $l\!\leq\! n$, that 
$$\forall t\!\in\!\omega^{n+1}\ 0^{n-l}\!\subseteq\! t\ \Rightarrow\ 
[t(0)0^{\alpha (0)}...t(n\! -\! 1)0^{\alpha (n-1)}t(n)0^p\!\subseteq\!\gamma\ \Rightarrow\ 
s_0t(0)...s_{n}t(n)\!\subseteq\! u(\gamma )].$$
We already proved it for $l\! =\! 0$. Assume that it is true for $l\! <\! n$, let $t\!\in\!\omega^{n+1}$ with 
$0^{n-l-1}\!\subseteq\! t$, and assume that 
$t(0)0^{\alpha (0)}...t(n\! -\! 1)0^{\alpha (n-1)}t(n)0^p\!\subseteq\!\gamma$. We set 
$\varepsilon\! :=\!\gamma\! -\! [t(0)0^{\alpha (0)}...t(n\! -\! 1)0^{\alpha (n-1)}t(n)0^p]$. Then by induction assumption on $l$ we get 
$$s_00...s_{n-l-1}0s_{n-l}t(n\! -\! l)...s_{n}t(n)\!\subseteq\! 
u[s^\omega_{\Sigma_{j<n-l} [1+\alpha (j)]}t(n\! -\! l)0^{\alpha (n-l)}...t(n\! -\! 1)0^{\alpha 
(n-1)}t(n)0^p\varepsilon ].$$
But by  induction assumption on $n$ we get, for each $i\!\in\!\omega$, 
$$s_00...s_{n-l-2}0s_{n-l-1}i\!\subseteq\! u[s^\omega_{\Sigma_{j<n-l-1} [1+\alpha (j)]}i0^{\alpha (n-l-1)}
t(n\! -\! l)0^{\alpha (n-l)}...t(n\! -\! 1)0^{\alpha (n-1)}t(n)0^p
\varepsilon ].$$
But $(u[s^\omega_{\Sigma_{j<n-l-1} [1+\alpha (j)]}i0^{\alpha (n-l-1)}
t(n\! -\! l)0^{\alpha (n-l)}...t(n\! -\! 1)0^{\alpha (n-1)}t(n)0^p
\varepsilon ])_{i\in\omega}\!\in\!\mathbb{A}_{\omega}$. This implies, for each $i\!\in\!\omega$, that  
$u[s^\omega_{\Sigma_{j<n-l-1} [1+\alpha (j)]}i0^{\alpha (n-l-1)}
t(n\! -\! l)0^{\alpha (n-l)}...t(n\! -\! 1)0^{\alpha (n-1)}t(n)0^p
\varepsilon ]$ begins with $s_00...s_{n-l-2}0s_{n-l-1}is_{n-l}t(n\! -\! l)...s_{n}t(n)$. In particular, this holds for $i\! =\! t(n\! -\! l\! -\! 1)$, and we are done.\bigskip

 It remains to choose $\alpha (n)\!\geq\! p$ such that 
$0^{\Sigma_{j\leq n} [1+\alpha (j)]}\! =\! s^\omega_{\Sigma_{j\leq n} [1+\alpha (j)]}$.

\vfill\eject

\noindent $\bullet$ If $s\!\in\!\omega^{\leq\omega}$, then we set 
$N[s]\! :=\!\hbox{\rm Card}\{ n\!\in\!\omega\mid s^\omega_n0\!\subseteq\! s\}$. Note that 
$N[\alpha ]\! =\!\omega$ if $\alpha\!\in\!\mathbb{G}$. By induction on $p$, we can construct $\beta (p)\!\in\!\omega$ such that $N[s_0\beta (0)...s_p\beta (p)s_{p+1}]\! =\! N[s_0]$. This implies that 
$N[s_0\beta (0)s_1\beta (1)...]\! =\! N[s_0]\! <\!\omega$, and 
$u[\beta (0)0^{\alpha (0)}\beta (1)0^{\alpha (1)}...]\!\notin\!\mathbb{G}$ by the previous point, which is absurd.\hfill{$\square$}

\section{$\!\!\!\!\!\!$ The proof in infinite dimension.}\indent

 Before proving Theorem 1.6, note first the following result:

\begin{thm} There is no $(X_0,\mathbb{A}_0)$, where $X_0$ is a metrizable compact space and 
$\mathbb{A}_0\!\in\!\ana (X_0^\omega )$, such that for any Polish space $X$, and for any 
$A\!\in\!\ana (X^\omega )$, exactly one of the following holds:\smallskip  

\noindent (a) $[X,A]\preceq_{B}[{\omega},\neg\Delta^\omega (\omega )]$.\smallskip  

\noindent (b) $[X_0,\mathbb{A}_{0}]\preceq_{c}[X,A]$.\end{thm}

\noindent\bf Proof.\rm\ Suppose towards a contradiction that such $(X_0,\mathbb{A}_0)$ exists. Note that $\mathbb{A}_0\!\not=\!\emptyset$, since otherwise we would have 
$[X_0,\mathbb{A}_{0}]\preceq_{B}[\omega ,\neg\Delta^\omega (\omega )]$. By Lemma 2.1, we now get some continuous $u\! :\! X_0\!\rightarrow\!\omega^\omega$ such that 
$\mathbb{A}_{0}\!\subseteq\! (u^\omega )^{-1}(\mathbb{A}_\omega )$. Then $u[X_0]$ will be a compact subset of $\omega^\omega$ and hence contained in some product 
$k_0\!\times\! k_1\!\times\! ...\!\subseteq\!\omega^\omega$, where the $k_i$'s are finite. Notice however that $(k_0\!\times\! k_1\!\times\! ...)^\omega\cap\mathbb{A}_\omega\! =\!\emptyset$, and thus 
$\mathbb{A}_{0}\!\subseteq\! (u^\omega )^{-1}[(k_0\!\times\! k_1\!\times\! ...)^\omega\cap
\mathbb{A}_\omega]\! =\!\emptyset$, which is a contradiction.\hfill{$\square$}\bigskip 

 Assume temporarily that there is a Polish space $X_0$ and $\mathbb{A}_0$ such that the end of the statement of Theorem 4.1 holds. By  Theorem 4.1, $X_0$ cannot be compact. Note that we may assume that $X_0$ is zero-dimensional, since there is a finer zero-dimensional Polish topology on $X_0$ (see 13.5 in [K]). This means that we can view $X_0$ as a closed subspace of $\omega^\omega$ (see 7.8 in [K]). As $X_0$ is not compact, the tree associated with this closed set (see 2.4 in [K]) is not finite 
splitting (see 4.11 in [K]). The proof of Theorem 1.6 will have the same scheme as the proof of Theorem 1.4. We have to build infinitely many $V_s$'s at some levels of the construction, since the tree associated with $X_0$ is not finite splitting. The only place where the proof of Theorem 1.4 does not work in infinite dimension is when we write ``$C\!\in\!\Ana (X)$".\bigskip
  
  The main modifications to make are the following:\bigskip

\noindent - As we have to build infinitely many $V_s$'s at some levels of the construction, it is not clear at all that $C$ remains $\Ana$, since $\Ana$ is not closed under infinite intersections. However, $\Ana$ is closed under $\forall^\omega$, and this will be enough. We will have to build the $V_s$'s uniformly in $s$ at each level of the construction to ensure that $C$ is $\Ana$, and it is possible. We will also ensure that there are only finitely many $U_{n,t}$'s at each level of the construction, to ensure that $C$ is $\Ana$.\bigskip

\noindent - The reason why Theorem 3 is true is that we cannot control all the diameters in 
$\mathbb{G}$ at each level of a construction that would give a map 
$u\! :\!\omega^\omega\!\rightarrow\!\mathbb{G}$. We will only control finitely many diameters, since we want $C$ to be $\Ana$. This is the reason why we will work in $\mathbb{G}$ instead of $\omega^\omega$. This gives the possibility to control only one diameter at each level of the construction among the $V_s$'s (and finitely many among the $U_{n,t}$'s). So the point in the proof of Theorem 1.6 is that we cannot build the $\Ana$ sets uniformly at each level of the construction and control all the diameters at the same time.\bigskip
 
\noindent\bf Proof of Theorem 1.6.\rm\ Note first that we cannot have (a) and (b) simultaneously, by Lemma 2.1.\bigskip

\noindent $\bullet$ Note that there is a recursive map $\tilde s\! :\!\omega\!\rightarrow\!\omega$ such that $\tilde s(l)$ codes $s^\omega_l$, i.e., $\tilde s(l)\! =\! I(s^\omega_l)$ (see the notation in the introduction). Indeed, there is a recursive map $\tilde\varphi\! :\!\omega\!\rightarrow\!\omega$ whose restriction to Seq is an increasing bijection from Seq onto $\omega$. Now $(\tilde\varphi\vert_{\hbox{\rm Seq}})^{-1}$ defines a recursive map $\tilde\psi_\omega\! :\!\omega\!\rightarrow\!\omega$. It remains to note that $\tilde s(l)\! =\! t$ is equivalent to
$$t\!\in\!\hbox{\rm Seq}\ \ \hbox{\rm and}\ \ \hbox{\rm lh}(t)\! =\! l\ \ \hbox{\rm and}\ \ 
\forall i\! <\! l\ [i\! <\!\hbox{\rm lh}[\tilde\psi_\omega (l)]\ \ \hbox{\rm and}\ (t)_i\! =\! (\tilde\psi_\omega (l))_i]
\ \ \hbox{\rm or}\ \ [i\!\geq\!\hbox{\rm lh}[\tilde\psi_\omega (l)]\ \ \hbox{\rm and}\ \ (t)_i\! =\! 0].$$
$\bullet$ We may assume that\bigskip

\noindent - The $X^{\omega^l}$'s are recursively presented Polish spaces, for $l\!\in\!\omega$.\bigskip

\noindent - The projections are recursive.\bigskip

\noindent - The maps $\Pi_l\! :\!\omega\!\times\! X^{\omega^l}\!\rightarrow\! X$ defined by 
$$\Pi_l[t,(x_s)_{s\in\omega^l}]\! =\! x\ \ \Leftrightarrow\ \ t\!\in\!\hbox{\rm Seq}\ \ 
\hbox{\rm and}\ \ \hbox{\rm lh}(t)\! =\! l\ \ \hbox{\rm and}\ \ x\! =\! x_{\overline{t}}$$ are partial recursive functions on $\{ t\!\in\!\omega\mid t\!\in\!\hbox{\rm Seq}\ \ \hbox{\rm and}\ \ \hbox{\rm lh}(t)\! =\! l\}\!\times\! X^{\omega^l}$, for $l\!\in\!\omega$.\bigskip

\noindent - The maps $\Pi'_l\! :\!\omega^2\!\times\! X^{\omega^l}\!\rightarrow\! X^\omega$ defined by 
$$\Pi'_l[n,t,(x_s)_{s\in\omega^l}]\! =\! (y_i)_{i\in\omega}\ \ \Leftrightarrow\ \ t\!\in\!\hbox{\rm Seq}\ \ 
\hbox{\rm and}\ \ n\! +\! 1\! +\!\hbox{\rm lh}(t)\! =\! l\ \ \hbox{\rm and}\ \ \forall i\!\in\!\omega\ \ y_i\! =\! x_{s^\omega_ni\overline{t}}$$ are partial recursive functions on 
$\{ (n,t)\!\in\!\omega^2\mid t\!\in\!\hbox{\rm Seq}\ \ \hbox{\rm and}\ \ n\! +\! 1\! +\!\hbox{\rm lh}(t)\! =\! l\}\!\times\! X^{\omega^l}$, for $l\!\in\!\omega$.\bigskip

\noindent - $A\!\in\!\Ana (X^\omega )$.\bigskip

\noindent $\bullet$ We set $\Phi\! :=\!\{ C\!\subseteq\! X\mid C\ \hbox{\rm is}\ A\hbox{\rm -discrete}\}$. As $\Phi$ is $\Ca$ on $\Ana$, the first reflection theorem ensures that if $C\!\in\!\Ana (X)$ is $A$-discrete, then there is $D\!\in\!\Borel (X)$ which is $A$-discrete and contains $C$.\bigskip

\noindent $\bullet$ By Lemma 2.3 we may assume that $U\!\not=\! X$, so that $Y\! :=\! X\!\setminus\! U$ is a nonempty $\Ana$ subset of $X$. The previous point gives the following key property:
$$\forall C\!\in\!\Ana (X)\ \ \ (\emptyset\!\not=\! C\!\subseteq\! Y\ \Rightarrow\ A\cap C^\omega \!\not=\!\emptyset ).$$
$\bullet$ We construct $(x_s)_{s\in\omega^{<\omega}}\!\subseteq\! Y$, 
$(V_s)_{s\in\omega^{<\omega}}\!\subseteq\!\Ana (X)$, and 
$(U_{n,t})_{(n,t)\in\omega\times\omega^{<\omega}}\!\subseteq\!\Ana (X^\omega )$ satisfying the following conditions:
$$\begin{array}{ll}
& (1)\ x_s\!\in\! V_s\!\subseteq\! Y\cap\Omega_X\ \ \hbox{\rm and}\ \ 
(x_{s_n^\omega it})_{i\in\omega}\!\in\! U_{n,t}\!\subseteq\! A\cap Y^\omega\cap
\Omega_{X^\omega}\hbox{\rm ,}\cr\cr
& (2)\ V_{sm}\!\subseteq\! V_{s}\ \ \hbox{\rm and}\ \ 
U_{n,tm}\!\subseteq\! U_{n,t}\hbox{\rm ,}\cr\cr
& (3)\ \hbox{\rm diam}_{d_X}(V_{s^\omega_l0})\!\leq\! 2^{-l}\ \ \hbox{\rm and}\ \ 
[s_n^\omega 0t\! =\! s^\omega_l0\ \Rightarrow\ 
\hbox{\rm diam}_{d_{X^\omega}}(U_{n,t})\!\leq\! 2^{-l}]\hbox{\rm ,}\cr\cr
& (4)\ \mbox{For any fixed $\vert s\vert$, the relation ``$x\!\in\! V_s$" is a $\Ana$ condition on 
$(x,s)$,}\cr\cr
& (5)\ \mbox{For any fixed $n$ and fixed $\vert t\vert$, the relation ``$(x_i)_{i\in\omega}\!\in\! U_{n,t}$" is a $\Ana$ condition on $[(x_i)_{i\in\omega},t]$.}
\end{array}$$
$\bullet$ Assume that this is done. Fix $\alpha\!\in\!\mathbb{G}$. Then 
$(V_{\alpha\vert p})_{p\in\omega}$ is a decreasing sequence of nonempty clopen subsets of 
$[\Omega_X ,{\it\Sigma}_X ]$ whose $d_X$-diameters tend to zero, so there is $u(\alpha )$ in their intersection. This defines $u\! :\!\mathbb{G}\!\rightarrow\! X$. Note that 
$d_X[x_{\alpha\vert p},u(\alpha )]\!\leq\!\hbox{\rm diam}_{d_X}(V_{\alpha\vert p})$, so that $u$ is continuous and $(x_{\alpha\vert p})_{p\in\omega}$ tends to $u(\alpha )$ in 
$[X ,{\it\Sigma}_X ]$.\bigskip

 If $(s^\omega_ni\gamma )_{i\in\omega}\!\in\!\mathbb{A}_\omega\cap\mathbb{G}^\omega$, then 
$(U_{n,\gamma\vert p})_{p\in\omega}$ is a decreasing sequence of nonempty clopen subsets of $[\Omega_{X^\omega} ,{\it\Sigma}_{X^\omega}]$ whose 
$d_{X^\omega}$-diameters tend to zero, so there is $(\alpha_i)_{i\in \omega}$ in their intersection. Note that $(\alpha_i)_{i\in \omega}\!\in\! A$. Moreover, the sequence $([x_{s^\omega_ni(\gamma\vert p)}]_{i\in \omega})_{p\in\omega}$ tends to $(\alpha_i)_{i\in \omega}$ in 
$[X^\omega ,{\it\Sigma}_{X^\omega}]$, and in 
$[X^\omega ,{\it\Sigma}_{X}^\omega]$ too. As $(x_{s^\omega_ni(\gamma\vert p)})_{p\in\omega}$ tends to $u(s^\omega_ni\gamma )$ in $[X ,{\it\Sigma}_X ]$, we get $u(s^\omega_ni\gamma )\! =\!\alpha_i$, for each $i\!\in\! \omega$. Thus $[u(s^\omega_ni\gamma )]_{i\in\omega}\!\in\! A$.\bigskip

\noindent $\bullet$ So it is enough to see that the construction is possible. If $V_\emptyset$ is any $\Ana$ set, then clearly (4) holds for $s$ of length $0$. Now suppose that $V_s$ has been defined for all $s\!\in\!\omega^{\leq l}$ and that (4) holds. Then in order to define $V_r$ for $r\!\in\!\omega^{l+1}$, while ensuring (4), we will let $V_{s^\omega_l0}\!\subseteq\! V_{s^\omega_l}$ be some chosen $\Ana$ set of diameter at most $2^{-l}$ (to be determined later on) and $V_{sm}\! :=\! V_s$ for all 
$sm\!\not=\! s^\omega_l0$. Then for $r\!\in\!\omega^{l+1}$ 
$$x\!\in\! V_r\Leftrightarrow (r\! =\! s^\omega_l0\ \mbox{ and }\ x\!\in\! V_{s^\omega_l0})\ \mbox{ or }\ 
(r\! =\! sm\!\not=\! s^\omega_l0\ \mbox{ and }\ x\!\in\! V_{s})\mbox{,}$$
which is $\Ana$ in $(x,r)$ by the induction hypothesis.\bigskip

 Similarly, if $U_{n,\emptyset}$ is any $\Ana$ set, then clearly (5) holds for $t$ of length $0$. Now suppose that $U_{n,t}$ has been defined for all $t\!\in\!\omega^{\leq k}$ and that (5) holds. Then in order to define $U_{n,r}$ for $r\!\in\!\omega^{k+1}$, while ensuring (5), we again split into two cases. If 
$s^\omega_n0r\! =\! s^\omega_n0t0\! =\! s^\omega_l0$, then $U_{n,r}\!\subseteq\! U_{n,t}$ will be some chosen $\Ana$ set of diameter at most $2^{-l}$ (to be determined later on). On the other hand, if 
$s^\omega_n0r\! =\! s^\omega_n0tm\!\not=\! s^\omega_l0$, then we set $U_{n,r}\! :=\! U_{n,t}$. Then for $r\!\in\!\omega^{k+1}$
$$(x_i)_{i\in\omega}\!\in\! U_{n,r}\Leftrightarrow\left\{\!\!\!\!\!\!\!\! 
\begin{array}{ll}
& \ (s^\omega_n0r\! =\! s^\omega_n0t0\! =\! s^\omega_l0\ \mbox{ and }\ 
(x_i)_{i\in\omega}\!\in\! U_{n,r})\cr
& \mbox{or}\cr 
& \ (s^\omega_n0r\! =\! s^\omega_n0tm\!\not=\! s^\omega_l0\ \mbox{ and }\ 
(x_i)_{i\in\omega}\!\in\! U_{n,t})\mbox{,}
\end{array}
\right.$$
which is $\Ana$ in $[(x_i)_{i\in\omega},r]$ by the induction hypothesis, since 
$s_n^\omega 0r\! =\! s^\omega_l0$ can hold for only finitely many 
$(n,r)\!\in\!\omega\!\times\!\omega^{<\omega}$.\bigskip

 Notice that in this way (2) and (3) are also satisfied, so it remains to define $V_{s^\omega_l0}$, 
$U_{n,\emptyset}$ and $U_{n,r}$ for $s^\omega_n0r\! =\! s^\omega_l0$ of diameter small enough such that (1) also holds.\bigskip

\noindent - As $Y$ is a nonempty $\Ana$ subset of $X$, we can choose $x_\emptyset\!\in\! Y\cap\Omega_X$, and set $V_\emptyset\! :=\! Y\cap\Omega_X$.\bigskip

\noindent - The key property applied to $V_\emptyset$ gives 
$(x_i)_{i\in\omega}\!\in\! A\cap V_\emptyset^\omega\cap\Omega_{X^\omega}$. We choose 
$U_{0,\emptyset}\!\in\!\Ana (X^\omega )$ such that 
$(x_i)_{i\in\omega}\!\in\! U_{0,\emptyset}\!\subseteq\! A\cap V_\emptyset^\omega\cap
\Omega_{X^\omega}$ and $\hbox{\rm diam}_{d_{X^\omega}}(U_{0,\emptyset})\!\leq\! 1$. Then we choose $V_0\!\in\!\Ana (X)$ such that $x_0\!\in\! V_0\!\subseteq\! V_\emptyset$ and 
$\hbox{\rm diam}_{d_{X}}(V_0)\!\leq\! 1$. Assume that $(x_s)_{\vert s\vert\leq l}$, 
$(V_s)_{\vert s\vert\leq l}$, and $(U_{n,t})_{n+1+\vert t\vert\leq l}$ satisfying (1)-(5) have been constructed, which is the case for $l\!\leq\! 1$.\bigskip

\noindent - We put\bigskip

\leftline{$C\! :=\!\Big\{\ x\!\in\! X\mid
\exists (y_s)_{s\in\omega^{l}}\!\in\! X^{\omega^l}\ \ y_{s^\omega_l}\! =\! x\ \ 
\hbox{\rm and}\ \ \forall s\!\in\!\omega^l\ \ y_s\!\in\! V_{s}\ \ \hbox{\rm and}\ \ \forall n\! <\! l\ \ 
\forall t\!\in\!\omega^{l-n-1}$}

\rightline{$(y_{s^\omega_{n}it})_{i\in\omega}\!\in\! U_{n,t}\ \Big\}.$}

\vfill\eject

 Then $x_{s^\omega_l}\!\in\! C$, by induction assumption. Moreover, $C\!\in\!\Ana$, by 
conditions (4) and (5) since $\Ana$ is closed under $\forall^\omega$. The key property applied to 
$C$ gives $(x_{s^\omega_li})_{i\in\omega}\!\in\! A\cap C^\omega\cap\Omega_{X^\omega}$. As 
$x_{s^\omega_lm}\!\in\! C$, there is 
$(x_{sm})_{s\in\omega^{l}\setminus \{ s^\omega_l\}}\!\subseteq\! X$ such that 
$x_{sm}\!\in\! V_{s}$ for each $s\!\in\!\omega^l$ and 
$(x_{s^\omega_{n}itm})_{i\in\omega}\!\in\! U_{n,t}$ for each $n\! <\! l$ and each 
$t\!\in\!\omega^{l-n-1}$. This defines $(x_s)_{s\in\omega^{l+1}}$.\bigskip

 We choose $U_{l,\emptyset}\!\in\!\Ana (X^\omega )$ such that $(x_{s^\omega_li})_{i\in\omega}\!\in\! U_{l,\emptyset}\subseteq\! A\cap V_{s^\omega_l}^\omega\cap\Omega_{X^\omega}$ and 
$\hbox{\rm diam}_{d_{X^\omega}}(U_{l,\emptyset})\!\leq\! 2^{-l}$, and $V_{s^\omega_l0}\!\in\!\Ana (X)$ such that $x_{s^\omega_l0}\!\in\! V_{s^\omega_l0}\!\subseteq\! V_{s^\omega_l}$ and 
$\hbox{\rm diam}_{d_{X}}(V_{s^\omega_l0})\!\leq\! 2^{-l}$. If 
$s_n^\omega 0r\! =\! s_n^\omega 0t0\! =\! s^\omega_l0$, then we choose 
$U_{n,r}\!\in\!\Ana (X^\omega )$ such that 
$\hbox{\rm diam}_{d_{X^\omega}}(U_{n,r})\!\leq\! 2^{-l}$ and 
$(x_{s_n^\omega ir})_{i\in\omega}\!\in\! U_{n,r}\!\subseteq\! U_{n,t}$.\hfill{$\square$}

\section{$\!\!\!\!\!\!$ The Baire-measurable natural extension in infinite dimension works.}\indent

 Theorem 1.7 is a consequence of Theorem 1.6, Lemma 2.1 and of the following result:

\begin{thm} $[\omega^\omega,\mathbb{A}_{\omega}]\preceq_{Bm}
[\mathbb{G},\mathbb{A}_{\omega}\cap\mathbb{G}^\omega ]$.\end{thm}

\noindent\bf Notation.\rm\ We define the following equivalence relation on the Baire space 
$\omega^\omega$, which is the analogous version of the usual equivalence relation $\mathbb{E}_0$ on the Cantor space $2^\omega$ (see [H-K-L]):
$$\alpha\ \mathbb{E}_0^{\omega^\omega}\beta\ \ \ \Leftrightarrow\ \ \ \exists m\!\in\!\omega\ \ \forall n\!\geq\! m\ \ ~\alpha (n)\! =\!\beta (n).$$

\begin{lem} There is a dense and $\mathbb{E}_0^{\omega^\omega}$-invariant $G_\delta$ subset 
$G$ of $\omega^\omega$ such that 
$$\forall\alpha\!\in\! G\ \ \forall l,m\!\in\!\omega\ \ \exists n\!\geq\! m\ \ \ s^\omega_nl\!\subseteq\!\alpha$$
(in particular, $G\!\subseteq\!\mathbb{G}$).\end{lem}

\noindent\bf Proof.\rm\ We set $G_0\! :=\!\{\alpha\!\in\!\omega^\omega\mid\forall l,m\!\in\!\omega\ \ 
\exists n\!\geq\! m\ \ \ s^\omega_nl\!\subseteq\!\alpha\}$. Note that $G_0$ is a dense $G_\delta$ subset 
of $\omega^\omega$. We also define, for $n,p\!\in\!\omega$, $f_n^p\! :\!\omega^\omega\!\rightarrow\!
\{\alpha\!\in\!\omega^\omega\mid\alpha (n)\! =\! p\}$ by 
$$f_n^p(\alpha )(m)\! :=\!\left\{\!\!\!\!\!\!\begin{array}{ll}
& \alpha (m)\ \ \mbox{if}\ \ m\!\not=\! n\mbox{,}\cr & \cr
& p\ \ \mbox{if}\ \ m\! =\! n.\end{array}\right.$$
Note that $f_n^p$ is onto, continuous, open, and has a clopen range. Then we set 
$$D\! :=\!\{ H\!\subseteq\!\omega^\omega\mid H\ \mbox{is a dense }G_\delta\}$$
and we define $\Phi\! :\! D\!\rightarrow\! D$ by $\Phi (H)\! :=\! H\cap\bigcap_{n,p\in\omega}\ (f_n^p)^{-1}(H)$. This allows us to define, for $q\!\in\!\omega$, $G_{q+1}\! :=\!\Phi (G_q)$, and we set 
$G\! :=\!\bigcap_{q\in\omega}\ G_q$. Note that $G\!\in\! D$. Moreover, if $\alpha\!\in\! G$ and 
$n,p\!\in\!\omega$, then $f_n^p(\alpha )\!\in\! G$. Indeed, let $q\!\in\!\omega$. Then 
$\alpha\!\in\! G_{q+1}\!\subseteq\! (f_n^p)^{-1}(G_q)$. Now if 
$\beta\ \mathbb{E}_0^{\omega^\omega}\alpha$, then there is $s\!\in\!\omega^{<\omega}$ such that 
$\beta\! =\! s(\alpha\! -\!\alpha\vert\vert s\vert )$ (which means that $s\!\subseteq\!\beta$ and $\alpha$, 
$\beta$ agree from the coordinate $\vert s\vert$ on). We set, for $i\!\leq\!\vert s\vert$, 
$\beta_i\! :=\! (s\vert i)(\alpha\! -\!\alpha\vert i)$, so that $\beta_0\! =\!\alpha$ and 
$\beta_{\vert s\vert}\! =\!\beta$. Note that $\beta_{i+1}\! =\! f_i^{s(i)}(\beta_i)$ for each $i\! <\!\vert s\vert$, by induction on $i$. This proves that $\beta_i\!\in\! G$ for each $i\!\leq\!\vert s\vert$, by induction on $i$. In particular, $\beta\!\in\! G$ which is $\mathbb{E}_0^{\omega^\omega}$-invariant. This finishes the proof since $G\!\subseteq\! G_0$.\hfill{$\square$}

\vfill\eject

\noindent\bf Notation.\rm\ For each $l\!\in\!\omega$, we define an oriented graph $G^\rightarrow_{l+1}$ on $\omega^{l+1}$ as follows:
$$s\ G^\rightarrow_{l+1}\ s'\ \ \Leftrightarrow\ \ \exists n\!\in\!\omega\ \ \exists i\!\not=\! 0\ \ 
\exists t\!\in\!\omega^{<\omega}\ \ (s,s')\! =\! (s^\omega_n0t,s^\omega_nit).$$
We denote by $G_{l+1}$ the symmetrization of $G^\rightarrow_{l+1}$.

\begin{lem} The graph $(\omega^{l+1},G_{l+1})$ is connected and acyclic.\end{lem}

\noindent\bf Proof.\rm\ We argue by induction on $l$. For $l\! =\! 0$, we have 
$$i\ G_{1}\ i'\ \Leftrightarrow\ (i\! =\! 0\ \hbox{\rm and}\ i'\!\not=\! 0)\ \hbox{\rm or}\ 
(i'\! =\! 0\ \hbox{\rm and}\ i\!\not=\! 0).$$ 
If $i\! <\! i'$, then $(i,0,i')$ is a $G_1$-walk from $i$ to $i'$ if $i\!\not=\! 0$, and 
$(i,i')$ is a $G_1$-walk from $i$ to $i'$ if $i\! =\! 0$. Thus $(\omega ,G_1)$ is connected. Now if 
$(i_j)_{j\leq L}$ is a $G_1$-cycle, then either $i_0\!\not=\! 0$ and $i_1\! =\! i_{L-1}\! =\! 0$, or 
$i_0\! =\! 0$ and $i_2\! =\! 0$. In both cases, this contradicts the fact that $(i_j)_{j\leq L}$ is a cycle. Thus 
$(\omega ,G_1)$ is acyclic.\bigskip

 Assume that the result is true for $l$. Note that 
$$si\ G_{l+2}\ s'i'\ \ \Leftrightarrow\ \ (s\! =\! s'\! =\! s^\omega_{l+1}\ \ \hbox{\rm and}\ \ i\ G_1\ i')\ \ 
\hbox{\rm or}\ \ (s\ G_{l+1}\ s'\ \ \hbox{\rm and}\ \ i\! =\! i').$$
We set, for $i\!\in\!\omega$, $E_i\! :=\!\{ ti\mid t\!\in\!\omega^{l+1}\}$. Note that $\omega^{l+2}$ is the disjoint union of the $E_i$'s, that the map $ti\!\mapsto\! t$ is an isomorphism from $(E_i,G_{l+2})$ onto 
$(\omega^{l+1},G_{l+1})$, and that the map $s^\omega_{l+1}i\!\mapsto\! i$ is an isomorphism from 
$(\{ s^\omega_{l+1}i\mid i\!\in\!\omega\},G_{l+2})$ onto $(\omega,G_{1})$. In particular, $(E_i,G_{l+2})$ is connected and acyclic, and $(\omega^{l+2},G_{l+2})$ is connected.\bigskip

 Now if $(t_j)_{j\leq L}$ is a $G_{l+2}$-cycle, then the sequence $[t_j(l\! +\! 1)]_{j\leq L}$ is not constant. There are $j_0\!\leq\! L$ minimal with $t_{j_0}(l\! +\! 1)\!\not=\! t_{0}(l\! +\! 1)$, and $j_1\! >\! j_0$ minimal with $t_{j_1}(l\! +\! 1)\! =\! t_{0}(l\! +\! 1)$. Note that 
$t_{j_0-1}\! =\! t_{j_1}\! =\! s^\omega_{l+1}t_0(l\! +\! 1)$. Thus $j_0\! =\! 1$ and $j_1\! =\! L$. 
If $t_{0}(l\! +\! 1)\!\not=\! 0$, then $t_{1}\! =\! t_{L-1}\! =\! s^\omega_{l+1}0$. If $t_{0}(l\! +\! 1)\! =\! 0$, then the sequence $[t_j(l\! +\! 1)]_{0<j<L}$ is constant and 
$t_{1}\! =\! t_{L-1}\! =\! s^\omega_{l+1}t_{1}(l\! +\! 1)$. In both cases, this contradicts the fact that $(t_j)_{j\leq L}$ is a cycle. Thus $(\omega^{l+2},G_{l+2})$ is acyclic.\hfill{$\square$}\bigskip

\noindent\bf Notation.\rm\ Lemma 5.3 and Theorem I.2.5 in [B] imply the existence, for each pair $\{ s,s'\}$ of distinct vertices in $\omega^{l+1}$, of a unique $s\! -\! s'$ path in $(\omega^{l+1},G_{l+1})$. We will call it $p^{l+1}_{s,s'}$. If $s\! =\! s'$, then we set $p^{l+1}_{s,s'}\! := <s>$. The proof of Lemma 5.3 shows that
$$p^{l+2}_{si,s'i'}\! =\!\left\{\!\!\!\!\!\!\!
\begin{array}{ll}
& <p^{l+1}_{s,s'}(0)i,...,p^{l+1}_{s,s'}(\vert p^{l+1}_{s,s'}\vert\! -\! 1)i>\hbox{\rm if}\ \ i\! =\! i'\hbox{\rm ,}\cr & \cr
& <p^{l+1}_{s,s^\omega_{l+1}}(0)i,...,
p^{l+1}_{s,s^\omega_{l+1}}(\vert p^{l+1}_{s,s^\omega_{l+1}}\vert\! -\! 1)i,s^\omega_{l+1}0,
p^{l+1}_{s^\omega_{l+1},s'}(0)i',...,
p^{l+1}_{s^\omega_{l+1},s'}(\vert p^{l+1}_{s^\omega_{l+1},s'}\vert\! -\! 1)i'>\cr 
& \ \ \ \ \ \ \ \ \ \ \ \ \ \ \ \ \ \ \ \ \ \ \ \ \ \ \ \ \ \ \ \ \ \ \ \ \ \ \ \ \ \ \ \ \ \ \ \ \ \ \ \ \ \ \ \ \ \ \ \ \ \ \ \ \ \ \ \ \ \ \ \ \ \ \ \ \ \ \ \ \ \ \ \ \ \ \ \ \ \ \ \ \ \ \ \ \ \ \ \ \ \ \ \ \ \ \ \ \ \ \ \ \ \ \ \hbox{\rm if}\ \ 0\!\not=\! i\!\not=\! i'\!\not=\! 0\hbox{\rm ,}\cr & \cr
& <p^{l+1}_{s,s^\omega_{l+1}}(0)i,...,
p^{l+1}_{s,s^\omega_{l+1}}(\vert p^{l+1}_{s,s^\omega_{l+1}}\vert\! -\! 1)i,
p^{l+1}_{s^\omega_{l+1},s'}(0)i',...,
p^{l+1}_{s^\omega_{l+1},s'}(\vert p^{l+1}_{s^\omega_{l+1},s'}\vert\! -\! 1)i'>\cr 
& \ \ \ \ \ \ \ \ \ \ \ \ \ \ \ \ \ \ \ \ \ \ \ \ \ \ \ \ \ \ \ \ \ \ \ \ \ \ \ \ \ \ \ \ \ \ \ \ \ \ \ \ \ \ \ \ \ \ \ \ \ \ \ \ \ \ \ \ \ \ \ \ \ \ \ \ \ \ \ \ \ \ \ \ \ \ \ \ \ \ \ \ \ \ \ \ \ \ \ \ \ \ \ \ \ \ \ \ \ \ \ \ \ \ \ \ \ \ \ \ \ \ \ \hbox{\rm otherwise.}
\end{array}\right.$$

\begin{lem} Let $\beta\!\in\!\omega^\omega$. Then 
$[[\beta]_{\mathbb{E}_0^{\omega^\omega}},\mathbb{A}_{\omega}\cap 
([\beta]_{\mathbb{E}_0^{\omega^\omega}})^\omega ]\preceq
[G,\mathbb{A}_{\omega}\cap G^\omega ]$.\end{lem}

\noindent\bf Proof.\rm\ We have seen that if $\alpha\ \mathbb{E}_0^{\omega^\omega}\beta$, then there is $s\!\in\!\omega^{<\omega}$ such that $\alpha\! =\! s(\beta\! -\!\beta\vert\vert s\vert )$. We will construct 
$u(\alpha )\!\in\! G$ by induction on $\vert s\vert$.\bigskip

\noindent $\bullet$ If $\vert s\vert\! =\! 0$, then we simply choose $u(\beta )\!\in\! G$.\bigskip

\noindent $\bullet$ If $\vert s\vert\! =\! 1$, then we choose $n_0\!\in\!\omega$ such that 
$s^\omega_{n_0}\beta (0)\!\subseteq\! u(\beta )$, and we set 
$$u[i(\beta\! -\!\beta\vert 1)]\! :=\! s^\omega_{n_0}i[u(\beta )\! -\! u(\beta )\vert (n_0\! +\! 1)]$$ 
if $i\!\not=\!\beta (0)$. Note that 
$u[i(\beta\! -\!\beta\vert 1)]\ \mathbb{E}_0^{\omega^\omega}u(\beta )\!\in\! G$, so that 
$u[i(\beta\! -\!\beta\vert 1)]\!\in\! G$. Moreover, we have 
$(u[i(\beta\! -\!\beta\vert 1)])_{i\in\omega}\!\in\!\mathbb{A}_\omega$.\bigskip

\noindent $\bullet$ Assume that $u(\alpha )\!\in\! G$ is constructed for $\vert s\vert\!\leq\! l\! +\! 1$, which is the case for $l\! =\! 0$. Let $\varphi\! :\!\omega^{l+1}\!\rightarrow\!\omega$ be a bijection with 
$\varphi (s^\omega_{l+1})\! =\! 0$.\bigskip

\noindent $\bullet$ We construct $(E_q)_{q\in\omega}\!\in\! [{\cal P}(\omega)]^\omega$ $\subseteq$-increasing such that 
$[\{\varphi^{-1}(p)\mid p\!\in\! E_q\},G_{l+1}]$ is connected for each $q\!\in\!\omega$ (see Lemma 5.3). We proceed by induction on $q$. We first set $E_0\! :=\!\{ 0\}$. Assume that $E_q$ is constructed.\bigskip

\noindent - If $E_q\! =\!\omega$, then we set $E_{q+1}\! :=\!\omega$.\bigskip
 
\noindent - If $E_q\!\not=\!\omega$, then we use the paths $p^{l+1}_{s,s'}$ defined after Lemma 5.3. We choose $r\!\in\!\omega\!\setminus\! E_q$ minimal for which there is $p\!\in\! E_q$ such that $\vert p^{l+1}_{\varphi^{-1}(p),\varphi^{-1}(r)}\vert\! =\! 2$. Such an $r$ exists since if $m\!\in\!\omega\!\setminus\! E_q$, then there is 
$i\! <\!\vert p^{l+1}_{s^\omega_{l+1},\varphi^{-1}(m)}\vert$ minimal such that 
$\varphi [p^{l+1}_{s^\omega_{l+1},\varphi^{-1}(m)}(i)]\!\notin\! E_q$, and 
$\vert p^{l+1}_{p^{l+1}_{s^\omega_{l+1},\varphi^{-1}(m)}(i\! -\! 1),
p^{l+1}_{s^\omega_{l+1},\varphi^{-1}(m)}(i)}\vert\! =\! 2$ since $i\! >\! 0$. As 
$[\{\varphi^{-1}(p)\mid p\!\in\! E_q\},G_{l+1}]$ is connected, and acyclic by Lemma 5.3, there is a unique 
$p\!\in\! E_q$ such that $\varphi^{-1}(p)\ G_{l+1}\ \varphi^{-1}(r)$. There are $n\!\leq\! l$, $i_0\!\not=\! 0$ and $t\!\in\!\omega^{l-n}$ such that $\{\varphi^{-1}(p),\varphi^{-1}(r)\}\! =\!\{ s^\omega_n0t,s^\omega_ni_0t\}$. We set 
$E_{q+1}\! :=\! E_q\cup\{\varphi (s^\omega_nit)\mid i\!\in\!\omega\}$.\bigskip

\noindent\bf Claim 1\it\ $\bigcup_{q\in\omega}\ E_q\! =\!\omega$.\rm\bigskip

 Indeed, let $r\!\in\!\omega\!\setminus\!\{ 0\}$. By induction on $k\!\in\!\omega$ we see that 
$\varphi [p^{l+1}_{s^\omega_{l+1},\varphi^{-1}(r)}(1+k)]\!\in\!\bigcup_{q\in\omega}\ E_q$. Thus $r$ is in 
$\bigcup_{q\in\omega}\ E_q$.\bigskip

 This allows us to define $q(s)\! :=\!\hbox{min}\{ q\!\in\!\omega\mid\varphi (s)\!\in\! E_q\}$, for 
$s\!\in\!\omega^{l+1}$.\bigskip 

\noindent\bf Claim 2\it\ Let $n\!\leq\! l$, and $t\!\in\!\omega^{l-n}$. Then there is 
$i\!\in\!\omega$ such that $q(s^\omega_nit)\! <\! q(s^\omega_njt)$ for each $j\!\not=\! i$. Moreover, $q(s^\omega_njt)\! =\! q(s^\omega_nj't)$ if $j,j'\!\not=\! i$.\rm\bigskip 

 Indeed, we argue by contradiction. Choose $i\!\not=\! j$ such that $q\! :=\! q^{n,t}_i\! =\! q^{n,t}_j$ is minimal among the $q^{n,t}_k$'s. By definition of $E_0$ we have $q\!\not=\! 0$. As 
 $\varphi (s^\omega_nit)\!\in\! E_q\!\setminus\! E_{q-1}$, we have $E_{q-1}\!\not=\!\omega$. This implies the existence of $n'\!\leq\! l$ and $t'\!\in\!\omega^{l-n'}$ such that 
$E_q\!\setminus\! E_{q-1}\!\subseteq\!\{\varphi (s^\omega_{n'}it')\mid i\!\in\!\omega\}$. Thus $s^\omega_nit$ and $s^\omega_njt$ differ at the coordinate $n'$, which implies that $n\! =\! n'$ and $t\! =\! t'$. By construction of $E_q$ there is $k\!\in\!\omega$ such that $s^\omega_nkt\!\in\! E_{q-1}$, which contradicts the minimality of $q$. This proves Claim 2.\bigskip

\noindent $\bullet$ We have to construct $u(sk[\beta\! -\!\beta\vert (l\! +\! 2)])\!\in\! G$ for 
$\vert s\vert\! =\! l\! +\! 1$ and $k\!\not=\!\beta (l\! +\! 1)$. We will construct 
$u(sk[\beta\! -\!\beta\vert (l\! +\! 2)])$ by induction on $q(s)$.\bigskip

\noindent - If $q(s)\! =\! 0$, then $s\! =\! s^\omega_{l+1}$ and we choose $n_1\!\in\!\omega$ such that 
$s^\omega_{n_1}\beta (l\! +\! 1)\!\subseteq\! u(s^\omega_{l+1}[\beta\! -\!\beta\vert (l\! +\! 1)])$, and we set 
$$u(s^\omega_{l+1}k[\beta\! -\!\beta\vert (l\! +\! 2)])\! :=\! s^\omega_{n_1}k
[u(s^\omega_{l+1}[\beta\! -\!\beta\vert (l\! +\! 1)])\! -\! u(s^\omega_{l+1}[\beta\! -\!\beta\vert (l\! +\! 1)])\vert (n_1\! +\! 1)]$$ 
if $k\!\not=\!\beta (l\! +\! 1)$. As before, $u(s^\omega_{l+1}k[\beta\! -\!\beta\vert (l\! +\! 2)])\!\in\! G$. 
Moreover, 
$[u(s^\omega_{l+1}i[\beta\! -\!\beta\vert (l\! +\! 2)])]_{i\in\omega}\!\in\!\mathbb{A}_\omega$.\bigskip

\noindent - Assume that $u(sk[\beta\! -\!\beta\vert (l\! +\! 2)])\!\in\! G$ is constructed for $q(s)\!\leq\! q$, which is the case for $q\! =\! 0$. If $q(s)\! =\! q\! +\! 1$, then $\varphi (s)\!\in\! E_{q+1}\!\setminus\! E_q$. This implies the existence of $n\!\leq\! l$, $t\!\in\!\omega^{l-n}$, $i_0\!\not=\! 0$ and of a unique 
$p\!\in\! E_q$ such that $\{\varphi^{-1}(p),s\}\! =\!\{ s^\omega_n0t,s^\omega_ni_0t\}$.\bigskip

 Note that $q[\varphi^{-1}(p)]\!\leq\! q$, so that 
$\beta_p\! :=\! u(\varphi^{-1}(p)k[\beta\! -\!\beta\vert (l\! +\! 2)])$ is defined and in $G$. We choose $n_{q+1}\!\in\!\omega$ such that 
$s^\omega_{n_{q+1}}[\varphi^{-1}(p)(n)]\!\subseteq\!\beta_p$, and we set 
$$u(s^\omega_{n}itk[\beta\! -\!\beta\vert (l\! +\! 2)])\! :=\! s^\omega_{n_{q+1}}i
[\beta_p\! -\!\beta_p\vert (n_{q+1}\! +\! 1)]$$ 
if $i\!\not=\!\varphi^{-1}(p)(n)$. This is licit by Claim 2, since only $\beta_p$ is defined among the 
$u(s^\omega_{n}itk[\beta\! -\!\beta\vert (l\! +\! 2)])$'s. As before, 
$u(s^\omega_{n}itk[\beta\! -\!\beta\vert (l\! +\! 2)])\!\in\! G$. Moreover, 
$[u(s^\omega_{n}itk[\beta\! -\!\beta\vert (l\! +\! 2)])]_{i\in\omega}\!\in\!\mathbb{A}_\omega$.\bigskip

\noindent $\bullet$ Now $u\! :\! [\beta]_{\mathbb{E}_0^{\omega^\omega}}\!\rightarrow\! G$ is constructed. Assume that $(s^\omega_ni\gamma )_{i\in\omega}\!\in\!\mathbb{A}_\omega\cap 
([\beta]_{\mathbb{E}_0^{\omega^\omega}})^\omega$. We can write $\gamma\! =\!\tilde t\delta$, where 
$\tilde t\!\in\!\omega^{<\omega}$, $\delta\! =\!\beta\! -\!\beta\vert (n\! +\! 1\! +\!\vert\tilde t\vert )$, and 
$\tilde t(\vert\tilde t\vert\! -\! 1)\!\not=\!\beta (n\! +\!\vert\tilde t\vert )$ if $\tilde t\!\not=\!\emptyset$. We have to check that $[u(s^\omega_ni\gamma )]_{i\in\omega}\!\in\!\mathbb{A}_\omega$. We may assume that 
$\tilde t\!\not=\!\emptyset$. We set $k\! :=\!\tilde t(\vert\tilde t\vert\! -\! 1)$ and also 
$t\! :=\!\tilde t\vert (\vert\tilde t\vert\! -\! 1)$. Then 
$(s^\omega_ni\gamma )_{i\in\omega}\! =\!(s^\omega_nitk\delta )_{i\in\omega}$, and Claim 2 provides $i$. Now the construction of $u$ shows that $[u(s^\omega_ni\gamma )]_{i\in\omega}\!\in\!\mathbb{A}_\omega$ (consider $l\! :=\! n\! +\!\vert t\vert$).\hfill{$\square$}\bigskip

\noindent\bf Proof of Theorem 5.1.\rm\ Using the axiom of choice, fix a selector $S\! :\!\omega^\omega\!\rightarrow\!\omega^\omega$ for $\mathbb{E}_0^{\omega^\omega}$, i.e., a map satisfying 
$\alpha\ \mathbb{E}_0^{\omega^\omega}\beta\ \ \Rightarrow\ \ S(\alpha )\! =\! S(\beta )\ \mathbb{E}_0^{\omega^\omega}\alpha$ for each $\alpha ,\beta\!\in\!\omega^\omega$ (see 12.15 in [K]). We can write
$$\omega^\omega\! =\! G\cup\bigcup_{\beta\in S[\omega^\omega ]\setminus G}\ 
[\beta ]_{\mathbb{E}_0^{\omega^\omega}}\mbox{,}$$
and this union is disjoint. By Lemma 5.4 there is $u_\beta\! :\! [\beta]_{\mathbb{E}_0^{\omega^\omega}}\!\rightarrow\! G$ such that 
$$\mathbb{A}_{\omega}\cap ([\beta]_{\mathbb{E}_0^{\omega^\omega}})^\omega
\!\subseteq\! (u_\beta^\omega)^{-1}(\mathbb{A}_{\omega}\cap G^\omega )\mbox{,}$$ 
for each $\beta\!\in\!\omega^\omega$. We define $u\! :\!\omega^\omega\!\rightarrow\!\mathbb{G}$ by 
$$u(\alpha )\! :=\!\left\{\!\!\!\!\!\!\!\begin{array}{ll}
& \alpha\ \ \mbox{if}\ \ \alpha\!\in\! G\mbox{,}\cr & \cr
& u_\beta (\alpha )\ \ \mbox{if}\ \ \alpha\!\in\! [\beta ]_{\mathbb{E}_0^{\omega^\omega}}\ \ \mbox{and}\ \ 
\beta\!\in\! S[\omega^\omega]\!\setminus\! G.
\end{array}\right.$$
Now let $U$ be an open subset of $\mathbb{G}$. Then 
$u^{-1}(U)\! =\! (G\cap U)\cup\bigcup_{\beta\in S[\omega^\omega ]\setminus G}\ u_\beta^{-1}(U)$. The set 
$G\cap U$ is a $G_\delta$ subset of $\omega^\omega$, and $\bigcup_{\beta\in S[\omega^\omega ]\setminus G}\ u_\beta^{-1}(U)\!\subseteq\!\omega^\omega\!\setminus\! G$ is meager. This proves that $u$ is Baire-measurable.\bigskip

 Now let $(s^\omega_ni\gamma )_{i\in\omega}\!\in\!\mathbb{A}_\omega$. Note that 
$s^\omega_ni\gamma\ \mathbb{E}_0^{\omega^\omega}s^\omega_nj\gamma$ if $i,j\!\in\!\omega$. If 
$s^\omega_n0\gamma\!\in\! G$, then 
$$[u(s^\omega_ni\gamma )]_{i\in\omega}\! =\! 
(s^\omega_ni\gamma )_{i\in\omega}\!\in\!\mathbb{A}_\omega\cap G^\omega\!\subseteq\!
\mathbb{A}_\omega\cap\mathbb{G}^\omega$$
since $G$ is $\mathbb{E}_0^{\omega^\omega}$-invariant. If $s^\omega_n0\gamma\!\notin\! G$, then there is $\beta\!\in\! S[\omega^\omega ]\!\setminus\! G$ such that 
$s^\omega_n0\gamma\!\in\! [\beta ]_{\mathbb{E}_0^{\omega^\omega}}$. In this case we have 
$(s^\omega_ni\gamma )_{i\in\omega}\!\in\!\mathbb{A}_\omega\cap ([\beta]_{\mathbb{E}_0^{\omega^\omega}})^\omega$. Thus $[u_\beta (s^\omega_ni\gamma )]_{i\in\omega}\!\in\!
\mathbb{A}_\omega\cap G^\omega\!\subseteq\!\mathbb{A}_\omega\cap\mathbb{G}^\omega$ and 
$[u(s^\omega_ni\gamma )]_{i\in\omega}\!\in\!\mathbb{A}_\omega\cap\mathbb{G}^\omega$. This finishes the proof.\hfill{$\square$}\bigskip

\noindent\bf Question.\rm\ Is it true that 
$[\omega^\omega,\mathbb{A}_{\omega}]\preceq_{B}[\mathbb{G},\mathbb{A}_{\omega}\cap\mathbb{G}^\omega ]$? This would imply that we can replace ``Baire measurable" with ``Borel" in Theorem 1.7.

\section{$\!\!\!\!\!\!$ References.}

\noindent [B]\ \ B. Bollob\'as,~\it Modern graph theory,~\rm 
Springer-Verlag, New York, 1998

\noindent [H-K-L]\ \ L. A. Harrington, A. S. Kechris and A. Louveau, A Glimm-Effros 
dichotomy for Borel equivalence relations,~\it J. Amer. Math. Soc.\rm ~3 (1990), 903-928

\noindent [K]\ \ A. S. Kechris,~\it Classical Descriptive Set Theory,~\rm 
Springer-Verlag, 1995

\noindent [K-S-T]\ \ A. S. Kechris, S. Solecki and S. Todor\v cevi\'c, Borel chromatic numbers,\ \it 
Adv. Math.\rm\ 141 (1999), 1-44

\noindent [M]\ \ Y. N. Moschovakis,~\it Descriptive set theory,~\rm North-Holland, 1980

\noindent [S]\ \ G. E. Sacks,~\it Higher Recursion Theory,~\rm 
Springer-Verlag, 1990\bigskip\bigskip\bigskip

\centerline{$\bullet$ Universit\' e Paris 6, Institut de Math\'ematiques de Jussieu, tour 46-0, bo\^\i te 186,}

\centerline{4, place Jussieu, 75 252 Paris Cedex 05, France.}

\centerline{dominique.lecomte@upmc.fr}\bigskip

\centerline{$\bullet$ Universit\'e de Picardie, I.U.T. de l'Oise, site de Creil,}

\centerline{13, all\'ee de la fa\"\i encerie, 60 107 Creil, France.}

\end{document}